\documentclass[12pt,leqno]{article}

\usepackage{amssymb}
\usepackage{amsmath}
\usepackage{color}
\usepackage{url}
\makeindex

\newcommand{\comment}[1]{}
\newcommand{\A}{\mbox{$\mathfrak A$}}
\newcommand{\B}{\mbox{$\mathfrak B$}}

\newcommand{\K}{\mbox{$\sf K$}}
\newcommand{\y}{\bar{y}}
\newcommand{\z}{\bar{z}}
\newcommand{\s}{\bar{s}}

\newtheorem{thm}{Theorem}
\newtheorem{lem}{Lemma}
\newtheorem{corollary}{Corollary}

\begin{document}

\title{Ultraproducts of continuous posets}
\author{Andr\'eka, H., Gyenis, Z.\ and N\'emeti, I.}%
\date{}
\maketitle

We dedicate this paper to Ervin Fried, teacher and friend. \bigskip

\begin{abstract} It is known that nontrivial ultraproducts of complete partially
ordered sets (posets) are almost never complete. We show that
complete additivity of functions is preserved in ultraproducts of
complete posets. Since failure of this property is clearly preserved
by ultraproducts, this implies that complete additivity of functions
is an elementary property.
\end{abstract}

An $n$-ary function $f$ in a poset is said to be {\it completely
additive} if $f(s_1,...,s_n)$ is the supremum of $\{ f(x_1,...,x_n)
: x_1\in X_1,...,x_n\in X_n\}$ whenever $s_1,...,s_n$ are the
suprema of $X_1,...,X_n$ respectively. Completely additive functions
are also called sup-preserving, or continuous.

\begin{thm}\label{main}
Assume that $\B_i$ are complete posets with completely additive
operations, for $i\in I$. The operations remain completely additive
in $\B=P\langle\B_i : i\in I\rangle\slash F$, the ultraproduct of
the $\B_i$'s modulo any ultrafilter $F$ on $I$.
\end{thm}

\noindent {\bf Proof.} First we consider the case of a unary
operation $f$, and then we will reduce the general case to this one.
Assume that $f_i$ is sup-preserving in $\B_i$ for all $i\in I$ and
let $f$ denote their ultraproduct. We note that a sup-preserving
function is always monotonic, thus each $f_i$ is monotonic, and so
$f$ is also monotonic.
\smallskip

\noindent \underline{Case I: $f$ is unary.}

Assume that $X\subseteq B$ and $s$ is the supremum of $X$ in $\B$.
We want to show that $f(s)$ is the supremum of $\{ f(x) : x\in X\}$
in $\B$. Let $f(X)$ denote $\{ f(x) : x\in X\}$. It is clear that
$f(s)$ is an upper bound of $f(X)$, since $f$ is monotonic. Assume
that $y$ is any upper bound for $f(X)$, we want to show that $y\ge
f(s)$.

Consider the formula $\alpha(x,s,y)$ defined as $x\le s\land f(x)\le
y$.  Let  $z=\sup A$ denote that $z$ is the supremum of $A$. We show
that
\begin{equation}\label{supa-eq}
s=\sup A, \text{ where }A = \{ x : \alpha(x,s,y)\}.
\end{equation}
Indeed, $s$ is an upper bound for $A$ since $\alpha(x,s,y)\to x\le
s$. Assume $z$ is any upper bound for $A$. By our assumptions that
$s=\sup X$ and $y$ is an upper bound for $f(X)$ we have $X\subseteq
A$. Thus $z$ is an upper bound for $X$ (since it is so for $A$),
hence $z\ge s$ since $s=\sup X$, as was desired.

Now we will use {\L}os Lemma, the Fundamental Theorem of
Ultraproducts, to show the existence of the analogous suprema in
almost all of the $\B_i$. Since $A$ is defined by the formula
$\alpha$, we get that \eqref{supa-eq} is expressible by a
first-order logic formula $\sigma(s,y)$, namely the following
formula will do:
\[ \forall z(\forall x(\alpha(x,s,y)\to x\le z)\to s\le z) .\]
(We omitted the part $\forall x(\alpha(x,s,y)\to x\le s)$ since this
follows directly from the definition of $\alpha$.) Let $\s\in s$,
$\y\in y$ be arbitrary (i.e., $s=\s/F$, $y=\y/F$), by the {\L}os
Lemma then we have $\{ i\in I : \B_i\models\sigma(\s_i,\y_i)\} \in F
$. Let $J$ denote this set and let $i\in J$. Let $A_i=\{x\in B_i :
\alpha(x,\s_i,\y_i)\}$, then $\B_i\models\sigma(\s_i,\y_i)$ implies
that $\s_i=\sup A_i$. By continuity of $f_i$ we get that
\begin{equation}\label{cont-eq}
f_i(\s_i) = \sup f_i(A_i) .
\end{equation}
This last statement can also be expressed by a first-order logic
formula, namely by the following $\varphi(s,y)$:
\[\forall z(\forall x(\alpha(x,s,y)\to f(x)\le z)\to f(s)\le z). \]
(We used $\alpha(x,s,y)\to f(x)\le y$.) Thus, by \eqref{cont-eq} we
have that $\B_i\models\varphi(\s_i,\z_i)$ for all $i\in J\in F$. By
the {\L}os Lemma we get $\B\models\varphi(s,y)$, i.e.,
\begin{equation}\label{su-eq}
f(s) = \sup f(A) .
\end{equation}
Since $y$ is an upper bound of $f(A)$ (by $\alpha(x,s,y)\to f(x)\le
y$), we get $f(s)\le y$ as was desired.\medskip\goodbreak

\noindent \underline{Case II: $f$ is $n$-ary for some $n\ge 2$.}

We reduce this case to Case I by using a straightforward
generalization of Theorem 1.6 (i$'$), (ii$'$) 
of \cite{JT51} that states that an operation is sup-preserving iff
it is sup-preserving in each of its coordinates. We write out the
simple proof because we are in a slightly different setting (in
\cite{JT51}, the poset is assumed to be a Boolean algebra).

Assume that $\A$ is any poset in which $f$ is an $n$-ary function.
For any elements $a_1,\dots,a_n\in A$ and $j\in\{ 1,\dots, n\}$ let
$f(a_1,\dots,a_{j-1},-,a_{j+1},\dots,a_n)$ denote the unary function
that takes any $z\in A$ to
$f(a_1,\dots,a_{j-1},z,a_{j+1},\dots,a_n)$. (When $j=1$ use
$f(-,a_{2},\dots,a_n)$ and when $j=n$ use $f(a_1,\dots,a_{n-1},-)$
in place of $f(a_1,\dots,a_{j-1},-,a_{j+1},\dots,a_n)$.) Call
$f(a_1,\dots,a_{j-1},-,a_{j+1},\dots,a_n)$ a unary instance of $f$.

\begin{lem}[{a version of \cite[Theorem 1.6]{JT51}}]\label{lem}
Assume that $\A$ is any poset in which $f$ is an $n$-ary operation.
The following are equivalent:
\begin{enumerate}
\item[\textnormal{(i)}]  $f$ is sup-preserving;
\item[\textnormal{(ii)}] all the unary instances of $f$ are sup-preserving.
\end{enumerate}
\end{lem}

\noindent {\bf Proof of Lemma~\ref{lem}.} Assume (i), then (ii)
holds by using $X_i=\{a_i\}$ when $i\ne j$ in the definition of $f$
being sup-preserving.

For simplicity and transparency, we write out the proof of
$\textnormal{(ii)} \Rightarrow \textnormal{(i)}$ for the case $n=2$.
Assume (ii) and let $X,Y\subseteq A$ with $s,z$ being the suprema of
$X,Y$ respectively. We have to show that $f(s,z)$ is the supremum of
$Z=\{ f(x,y) : x\in X, y\in Y\}$. By our assumption (ii) we have
that $f(x,z)$ is the supremum of $Z_x=\{ f(x,y) : y\in Y\}$, for all
$x\in X$. Also, $f(s,z)$ is the supremum of $V=\{ f(x,z) : x\in
X\}$. Thus, it is enough to show that the supremum of $V$ is the
same as the supremum of $Z$. Clearly, any upper bound $b$ of $V$ is
an upper bound of $Z$, because $f(x,y)\le f(x,z)\le f(s,z)\le b$.
Assume that $b$ is an upper bound of $Z$, we want to show that it is
an upper bound of $V$. Indeed,  $f(x,z)$ is the supremum of
$Z_x\subseteq Z$, so $f(x,z)\le b$ since $b$ is an upper bound of
$Z\supseteq Z_x$, and we are done with proving the equivalence of
(i) and (ii). Lemma\ref{lem} has been proved.\hfill$\Box$ \medskip

Now we can complete the proof of Theorem~\ref{main} itself. Assume
that $f_i$ is a sup-preserving operation in each $\B_i$. We want to
show that all the unary instances of their ultraproduct $f$ are
sup-preserving in $\B$. For transparency, we write out the proof for
$n=2$. Let $y\in B$ and $\y\in y$. Then $f(y,-)$ is the ultraproduct
of the unary functions $f_i(\y_i,-)$ of $\B_i$ (for $i\in I$). All
these are sup-preserving by Lemma~\ref{lem} and our assumption that
$f_i$ of $\B_i$ are sup-preserving. By the already proven Case I
then we get that their ultraproduct, $f(y,-)$, is sup-preserving,
too. Similarly for $f(-,y)$ for any $y\in B$. By Lemma~\ref{lem}
then $f$ is sup-preserving in $\B$.  This completes the proof of
Theorem~\ref{main}.\hfill$\Box$ \medskip

\begin{corollary}\label{cor}
Complete additivity of a function can be expressed with a
first-order logic formula.
\end{corollary}

\noindent {\bf Proof of Corollary~\ref{cor}.} Let $\K$ denote the
class of all posets with an $n$-place continuous function $f$.
Theorem~\ref{main} states that $\K$ is closed under ultraproducts.
Now, the complement of $\K$ consists of those posets where $f$ is
not continuous. Failure of continuity of $f$ is expressible with an
existential second-order logic formula. Namely $\exists X\exists
c(c=\sup X\land f(c)\ne\sup f(X))$ expresses that the unary $f$ is
not continuous, and a similar formula can be stated for $n$-place
functions $f$. It is known, and can easily be seen, that these
second-order formulas are preserved by ultraproducts. So, the
complement of $\K$ is also closed under ultraproducts. This means
that $\K$ is finitely axiomatizable, i.e., there is a first-order
logic formula $\varphi$ such that $\K$ consists exactly of those
structures in which $\varphi$ is valid. This $\varphi$ expresses
complete additivity of the function $f$. \hfill$\Box$

\bigskip

We close the paper with some remarks.

(a) Our original paper contained only Theorem~\ref{main}. Then we
saw a theorem due to Holliday and Litak, \cite[Lemma 5.13]{H15},
stating a first-order formula that is equivalent with continuity of
operators in Boolean algebras. This made us realize that our
Theorem~\ref{main} easily implies a generalization of their theorem
to arbitrary posets, see Corollary~\ref{cor}. This corollary also
answers a question of the referee we were thinking about. Holliday
and Litak got their result about the same time as we got our
Theorem~\ref{main}, see \cite{HL15}. They use the formula
\[\forall xy[x\cdot f(y)\ne 0\to\exists z(0<z\le y\land\forall y'(0<y'\le z\to x\cdot f(y')\ne 0))]\]
to express continuity of $f$ in Boolean algebras. This formula
cannot be generalized to arbitrary posets. We get, instead, the
following formula from the proof of Theorem~\ref{main} to express
continuity in arbitrary posets:
\[ \forall xy[(y=\sup\{ z<y : f(z)\le x\}\to f(y)\le x)\land(x\le y\to f(x)\le f(y))] .\]

(b) Let us call a function \emph{quasi-complete} if it preserves the
suprema of non\-empty sets. Thus, quasi-complete functions are not
required to take the smallest element, which when it exists is the
supremum of the empty set, to itself, in other words, they are not
necessarily normal. Givant pointed out to us that the proof of
Theorem~\ref{main} goes through for quasi-complete operations, with
the following modification of the proof. Since we want to show
quasi-completeness of $f$, we assume $X\ne\emptyset$. Then
$A\ne\emptyset$ by $X\subseteq A$, and hence $A_i\ne\emptyset$ for
almost all $i\in I$. So we can use quasi-completeness of $f_i$ to
infer $f(\s_i)=\sup f(A_i)$ from $\s_i=\sup A_i$, and this is the
only place in the proof of Theorem~\ref{main} where we used
completeness of $f_i$. On the other hand, the property of an
operator being normal is obviously closed under ultraproducts. As
Givant has observed (in the context of Boolean algebras with
operators), an operator is complete if and only if it is
quasi-complete and normal. Consequently, Theorem~\ref{main} follows
at once from the version of Theorem~\ref{main} that is formulated
for quasi-complete operators.

(c) Our theorem has an application in considering ultraproducts of
complex algebras of relational structures. Let us note that complex
algebras are always complete Boolean algebras with sup-preserving
operations. Suppose one starts with a system of relational
structures, then forms an ultraproduct of their complex algebras.
Now, this ultraproduct most likely is not complete. Givant has shown
in \cite[Theorem 1.35, pp.\ 56--60]{G14} 
that the completion of this ultraproduct is canonically isomorphic
to the complex algebra of the corresponding ultraproduct of the
original relational structures. He noticed that for his proof to be
valid, one must know that the ultraproduct of the complex algebras
is a Boolean algebra with quasi-complete operators in order to be
able to form the completion of this algebra. We proved the present
Theorem~\ref{main} at his request, and this makes Givant's proof
complete.

\bigskip\bigskip\bigskip

\noindent
Alfr\'ed R\'enyi Institute of Mathematics, Hungarian
Academy of Sciences\\
Budapest, Re\'altanoda st.\ 13-15, H-1053 Hungary\\
andreka.hajnal@renyi.mta.hu,\\ gyenis.zalan@renyi.mta.hu,\\
nemeti.istvan@renyi.mta.hu

\end{document}